\documentclass[11pt,reqno]{amsart}

\usepackage{amssymb,latexsym,verbatim,pdfsync,color,mathrsfs,graphicx,caption,epstopdf,subfigure} 

\newcommand{\bbC}{{\mathbb C}}
\newcommand{\bbD}{{\mathbb D}}
\newcommand{\bbN}{{\mathbb N}}
\newcommand{\bbR}{{\mathbb R}}
\newcommand{\bbT}{{\mathbb T}}
\newcommand{\bbX}{{\mathbb X}} 
\newcommand{\bbZ}{{\mathbb Z}} 
\newcommand{\bbH}{{\mathbb H}} 

\def\bZ{{\mathbf Z}}

\def\cA{{\mathcal A}}

\def\cC{{\mathcal C}}
\def\cD{{\mathcal D}}

\def\cF{{\mathcal F}}

\def\cH{{\mathcal H}}

\def\cM{{\mathcal M}}
\def\cN{{\mathcal N}}

\def\cP{{\mathcal P}}

\def\cR{{\mathcal R}}
\def\cS{{\mathcal S}}

\def\cW{{\mathcal W}}

\def\cF{{\mathscr F}}

\def\sS{{\mathscr S}}
\def\sT{{\mathscr T}}

\def\Re{\operatorname{Re}}
\def\Im{\operatorname{Im}}

\def\la{\langle}
\def\ra{\rangle}

\def\eps{\varepsilon}
\def\z{\zeta} 
\def\vp{\varphi}
\def\ov{\overline}
\def\p{\partial}

\def\ms{\medskip}

\def\MS{M\"untz--Sz\'asz}

\def\pt{{\textstyle \frac{\pi}{2} }}

\def\Hol{\operatorname{Hol}}

\def\Log{\operatorname{Log}}

\newtheorem{thm}{Theorem}[section]
\newtheorem{prop}[thm]{Proposition}
\newtheorem{cor}[thm]{Corollary}

\newtheorem{remarks}[thm]{Remarks}

\begin{document}

\title[(Ir-)Regularity of canonical projection
operators]{(Ir-)Regularity of canonical projection operators on some 
weakly pseudoconvex domains}
\author[A. Monguzzi, M. M. Peloso]{
Alessandro Monguzzi, Marco M. Peloso}
\address{Dipartimento di Matematica, Universit\`a  di
  Milano-Bicocca, Via. R. Cozzi 55, 20126 Milano,  Italy}
\address{Dipartimento di Matematica, Universit\`a degli Studi di
  Milano, Via C. Saldini 50, 20133 Milano, Italy}
\email{{\tt alessandro.monguzzi@unimib.it}}
\email{{\tt marco.peloso@unimi.it}}
\thanks{{\em Math Subject Classification} 32A25, 32A36, 30H20.}

\begin{abstract}

In this paper we discuss some recent results concerning the regularity
and irregularity of the Bergman and Szeg\H o projections on some
weakly pseudoconvex domains that have the common feature to possess a
nontrivial Nebenh\"ulle.
\end{abstract}

\maketitle


\section*{Introduction}
\label{sec1}
In this note we survey some recent results on the analysis of canonical
projection operators, such as the Bergman and Szeg\H o projections,
on a family of domains that present some
pathological behavior.  These domains  have the common feature to
possess
a nontrivial {\em Nebenh\"ulle}, and they essentially are 
 the worm domain of K. Diederich and
 J.E. Forn\ae ss, the Hartogs triangle and some of its variants, and 
 some model worm domains introduced by C. Kiselman and studied, among
others by D. Barrett, S. Krantz and the authors of this note.

This note is an extended version of a seminar given by the second
named author at the Dipartimento di Matematica dell'Universit\`a della
Basilicata.  He wishes to thank such department and in particular
E. Barletta and
S. Dragomir for the kind invitation and the great hospitality.
\ms

The worm domain $\mathcal{W}_\mu$ was   introduced by 
K. Diederich and J.E. Forn\ae ss in \cite{DF}.
$$
\cW_\mu=\left\{(z_1,z_2)\in\bbC^2: |z_1-e^{i\log
    |z_2|^2}|^2<1-\eta(\log|z_2|^2)\right\} 
$$
where 
$\eta$ is smooth, even, convex, 
vanishing on $[-\mu, \mu]$, with $\eta(a)=1$, and $\eta'(a)>0$.  These
properties of $\eta$ imply that
$\cW_\mu$ is smooth, bounded and pseudoconvex.  Morevover $\cW_\mu$ is 
strictly pseudoconvex
at all points $(z_1,z_2) \in \partial\cW_\mu$ with $z_1\neq
0$. 
The set of points on the boundary $\cA= \{ (0,z_2):\, \big|\log |z_2|^2\big|\le
\mu\}$ is the {\em critical annulus}.

In \cite{DF} the following important features of $\cW_\mu$ were shown: 
\begin{itemize}
\item[(I)] $\cW_\mu$ has non-trivial Nebenh\"ulle (that is, there exists no
  neighborhood basis of pseudoconvex domains for $\cW_\mu$) [Diederich-Forn\ae
  ss];
  \item[(II)] 
 $\cW_\mu$  does not admit any plurisubharmonic defining function
 (that is,  a defining
function that is plurisubharmonic on the boundary).
\end{itemize}

Concerning ${\rm (I)}$, 
by the Hartogs's extension phenomenon indeed, it follows  that if $f$ is holomorphic
in a neighborhood of 
$$
\big\{ (0,z_2): \big| \log|z_2|^2\big|\le \pi
\big\} \cup 
\big\{ (z_1,z_2): \big| \log|z_2|^2\big|= \pi,\ 
\big| z_1 -e^{i\mu}\big| \le1 
\big\} \, ,
$$
then $f$  is also holomorphic in a neighborhood of the set 
$$
\big\{ (z_1,z_2): \big| \log|z_2|^2\big|\le \pi,\ 
\big| z_1 -e^{i\mu}\big| \le1 
\big\} \,.
$$
Regarding ${\rm (I)}$
we recall that, given a domain $\Omega\subseteq\bbC^n$,
a continuous function $\vp:\Omega\to (-\infty,0)$ is called 
a {\em bounded exhaustion function} if for all $c<0\in\bbR$, 
$$
\ov{\vp^{-1} (-\infty,c) }\cap \p\Omega =\emptyset \,.
$$ 
In \cite{DF}
Diederich and Forn\ae ss  proved that if
$\Omega\subseteq\bbC^n $ is smooth, bounded and pseudoconvex 
with defining function $\rho$, then 
 there exists $\tau\in(0,1]$ such that  $-(-\rho)^\tau$ is a
bounded stricly plurisubharmonic exhaustion function. 
Such an exponent $\tau=\tau_\rho$ is called a $\operatorname{DF}${\em -exponent} for
the defining function $\rho$.  We set
$$
\operatorname{DF}(\Omega) 
=\sup \big\{\tau_\rho:\, \rho\text{ defining function of }\Omega \big\}  \,,
$$
 and we call this value the {\em Diederich--Forn\ae ss index} of
 $\Omega$.
 In \cite{DF}, Diedirich and Forn\ae ss  proved that 
 $$
\operatorname{DF}(\cW_\mu)\le\pi/2\mu\,.
$$

\begin{figure}
\begin{center}
\includegraphics[width=7.4cm]{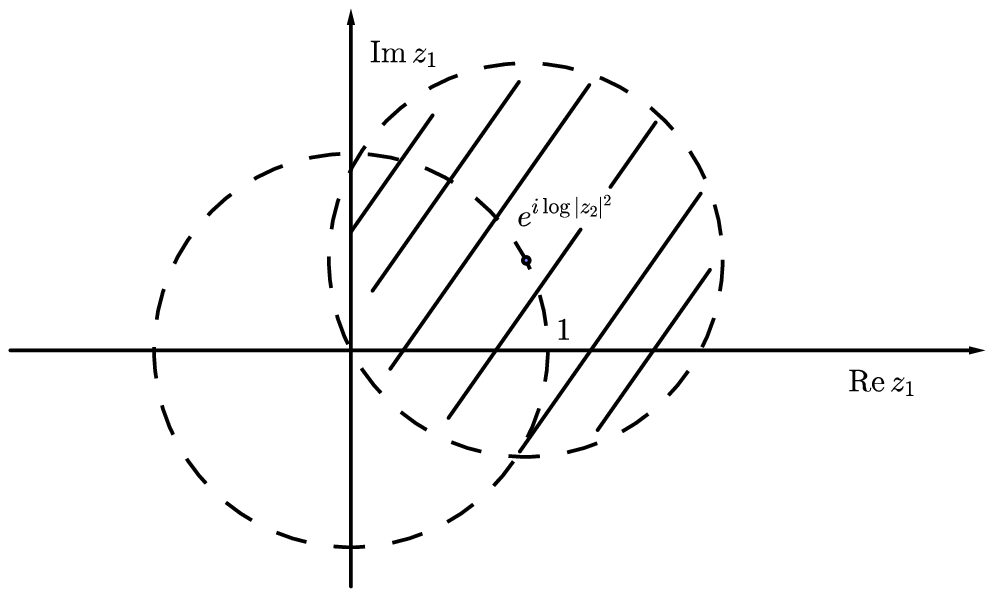}
\end{center}
\vspace{-0.4cm}
\caption{\footnotesize{Representation of $\cW_\mu$ in the $(\Re z_1, \Im z_1)$-plane.}}
\end{figure}

Since its appearance, research on the properties of  the worm domains
remained dormant for a number of years.   We now consider a still open
fundamental problem in analysis and geometry of several complex
variables:  Given  $\cD_1,\cD_2$ bounded, smooth,
pseudoconvex domains and  a biholomorphic mapping
$\Phi:\cD_1\to \cD_2$, does $\Phi$ extend smoothly to a diffeormophism of the
boundaries?
We denote by $\p\cD$ the topological boundary of a given domain
$\cD$.

Given a domain $\Omega$, let $A^2(\Omega)
=L^2(\Omega)\cap \operatorname{Hol}(\Omega)$
be the Bergman space, and 
$P_\Omega:L^2(\Omega)\to A^2(\Omega)$ be the Bergman projection. 
A celebrated theorem by S. Bell and E. Ligocka \cite{Bell-Ligocka},
and later improved by S. Bell \cite{Bell},  says that
this is the case if one of the two domains satisfies (say $\cD_1$) the so-called
Condition (R):
\begin{equation*}
P_{\cD_1} : C^\infty(\overline{\cD_1})\to  C^\infty(\overline{\cD_1})
\qquad
\text{is bounded}\,. \tag{R}
\end{equation*}

In \cite{Ba-Acta}  D. Barrett showed that, writing $\cW$ in place of
$\cW_\mu$ for short and denoting the Sobolev space on $\cW$ by
$W^{2,s}(\cW)$, 
$$ 
P_{\cW} :W^{2,s} (\cW)\not\to W^{2,s}(\cW)
$$
if $s\ge\pi/2\mu$.   Hence, $P_{\cW}$ fails to preserve the
$L^2$-Sobolev space $W^{2,s}(\cW)$, when $s$ is greater or equal to the
reciprocal of the windings of the domain $\cW$.

Although Barrett's result constituted a major breakthrough, it  did
not imply that $\cW$
failed to satisfy
Condition (R).  It was M. Christ in \cite{Christ} to prove that  the
Neumann operator $\cN$
on $\cW$ does not preserve
 $C^\infty(\overline{\cW})$; hence 
$\cW$ does not satisfy
Condition (R). In fact, by  a theorem of Boas--Straube
\cite{BS-equivalence}, on any given smoothly bounded pseudoconvex
domain $\Omega$, 
$\cN$ is globally (exactly)
  regular if and only if $P_\Omega$ is globally (exactly)
  regular.
  We say that $P_\Omega$ is exactly regular if
  $P_\Omega:W^{2,s}(\Omega)\to W^{2,s}(\Omega)$ is bounded for all $s>0$. We say
that $P_\Omega$ is globally regular if given any $s>0$, there exists
$q=q(s)$ such that
  $P_\Omega:W^{2,s+q(s)}(\Omega)\to W^{2,s}(\Omega)$ is bounded.  In
  particular, if $P_\Omega$ is globally regular, then $P_{\Omega} :
  C^\infty(\overline{\Omega})\to  C^\infty(\overline{\Omega})$ is
  bounded, that is, $\Omega$ satisfies condition (R).

The problem of the regularity of the Bergman projection on worm
domains 
has been object of active and intense research and 
 we mention in particular
\cite{KP-Houston, MR2904008,  MR3424478, CS, KPS1, KPS2}.    We also refer the reader to \cite{MR2393268} for a detailed
account on the subject. 
\ms

Main goal of this note is to report on recent progress on the analysis
of the  (ir-)regularity of  the
boundary analogue of the Bergman 
projection, that is, the Szeg\H{o} projection.
Given a smoothly bounded domain $\Omega=\{ z:\rho(z)<0\}
\subseteq\bbC^n$, the Hardy space $H^2(\Omega,d\sigma)$ is 
defined as 
$$
H^2(\Omega,d\sigma) =\big\{ f\in\Hol(\Omega): \sup_{\eps>0}
\int_{\p\Omega_\eps} |f|^2 d\sigma_\eps <\infty
\big\} \,,
$$
where $\Omega_\eps=\{ z:\rho(z)<-\eps\}$ and $d\sigma_\eps$ is the
induced surface measure on $\p\Omega_\eps$. 

Then, $H^2(\Omega,d\sigma)$ can be 
identified with a closed subspace of $L^2(\p\Omega,d\sigma)$, that we
denote by $H^2(\p\Omega,d\sigma)$, where $\sigma$ is the induced
surface measure on $\p\Omega$.  The
Szeg\H o projection is the orthogonal projection 
$$
S_\Omega : L^2(\p\Omega,d\sigma) \to H^2(\p\Omega,d\sigma) \,;
$$
see \cite{St2}.   \ms

The note is organized as follows. In Section \ref{2} we recall some
further  noticeable results concerning the worm domain, some of its
generalizations and some ideas involved in the proofs of such
results.  In Section \ref{Har-sec} we discuss the case of Szeg\H o
projections on worm domains.   Section \ref{Hartogs-sec} is devoted
to the case of another class of domains, the so-called Hartogs
triangles.  In Section \ref{MS-sec} we present an interesting problem
in the theory of $1$-dimensional Bergman spaces that arose in the
study of orthogonal sets in the Bergman space of the truncated worm
domain. We conclude this report with some final remarks and open
questions. 
\ms

\section{Generalizations and
  some open problems on
  the worm
  domain} \label{2}

The worm domain $\cW$ is still up to today the only known example of a
smoothly bounded pseudoconvex domain on which Condition (R) fails.
Thus, it is a natural testing ground for the validity of the
extendebility to diffeomorphism to the boundary of biholomorphic
mappings.   The first class of mappings that one is naturally led to
consider are the biholomorphic self-maps of $\cW$, that is, the {\em
  automorphisms} of $\cW$, $\operatorname{Aut}(\cW)$.
Clearly, the maps $\Phi(z_1,z_2)=(z_1,e^{i\theta}z_2)$ are in
$\operatorname{Aut}(\cW)$ and extend smoothly to the boundary.  The
obvious question is: Are
  there any others?   In \cite{SoChinChen}, the author studied
the automorphisms group
$\operatorname{Aut}(\cW)$, and claimed that 
 this is the case. Unfortunately, it is generally accepted that
 there is a gap  in the proof and
it has not been fixed.   Thus, the very interesting question of
characterizing the automorphism group 
$\operatorname{Aut}(\cW)$ is an open and fundamental question.
\ms

Before going any further, we point out that in \cite{MR2904008} D. Barrett and
S. \c{S}ahuto{\u{g}}lu constructed a higher dimensional analogue of
the worm domain.  Let $n\ge3$ and for $z\in\bbC^n$ we write
$z=(z_1,z',z_n) \in\bbC\times\bbC^{n-2}\times\bbC$.  For
$\lambda,\mu>0$ define
\begin{equation}\label{worm-3d}
  \cW_{\lambda,\mu} = \big\{ (z_1,z',z_n)
  \in\bbC\times\bbC^{n-2}\times\bbC: \,
  \big| z_1-e^{i\lambda \log|z_n|^2} \big|^2 < 1- |z'|^2
  -\widetilde\eta(\log|z_n|^2) 
  \big\} \,,
\end{equation}
where $\widetilde\eta$ is a {\em  particular}, explicit, smooth function which is
identically
$0$ when $e^{-1/2}\le |z_n| \le e^{\mu/2}$.  The function
$\widetilde\eta$ is chosen in such a way the domain is smoothly
bounded and pseudoconvex, and it is strongly pseudoconvex except at
the critical annulus
\begin{multline*}
\cA = \big\{ (z_1,z',z_n) \in\bbC\times\bbC^{n-2}\times\bbC\cap
\p\cW_{\lambda,\mu} :\, z_1=0, z'=0 \big\} \\
= \big\{ (0,0,z_n) \in\bbC\times\bbC^{n-2}\times\bbC :
e^{-1/2}\le |z_n| \le e^{\mu/2}\big\}
\,.
\end{multline*}

Barrett and
S. \c{S}ahuto{\u{g}}lu proved that the Bergman projection
$P_{\cW_{\lambda,\mu}}$ fails to preserve the Sobolev spaces
$W^{p,s}$, with $p\in[1,\infty)$ and $s\ge0$, hence
including the cases $p\neq2$, when
$$
s\ge \frac{\pi}{2\lambda\mu} +n\Big(\frac12-\frac1p\Big) \,.
$$

What is extremely interesting to notice here is that the Bergman
projection becomes irregular if either the winding is too ``long''
(i.e. when $\mu$ is large), or is too ``fast'' (i.e. when $\lambda$ is
large). \ms

For simplicity of presentation, we restricted ourselves to the
$2$-dimensional case, that is, to the domain $\cW=\cW_\mu$.  However, we point
out that the discussion that follows is also valid for the higher
dimensional cases of the domains $ \cW_{\lambda,\mu}$.

Instrumental to Barrett's proof of the irregularity of $P_\cW$ were
two unbounded model worm domains, that we denote by
$D_\mu$ and $D'_\mu$, where
\begin{equation*} 
D_\mu =
\Big\{ (z_1,z_2)\in\bbC^2:\, \Re\big(z_1e^{-i\log|z_2|^2}\big)>0,\ 
\big| \log|z_2|^2\big|<\mu  \Big\}, \quad \mu>\pi,
\end{equation*}
and
$$ 
D'_\mu=\Big\{(z_1,z_2)\in\bbC^2: \big|\Im z_1-\log|z_2|^2\big|<\frac\pi2,
\big|\log|z_2|^2\big|<\mu \Big\} \,.
$$

\begin{figure}[h]
\begin{center}
\includegraphics[width=8.7cm]{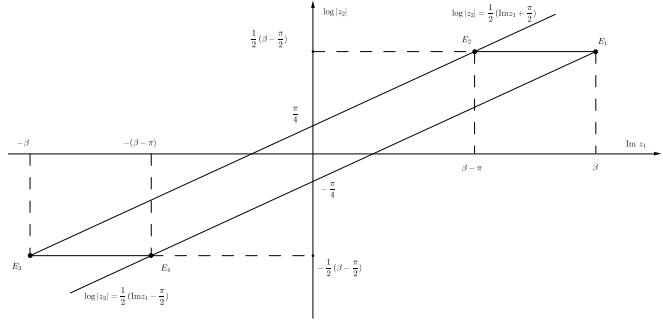}
\end{center}
\caption{\footnotesize{Representation of $D'_\mu$ in the $(\Im z_1, \log|z_2|)$-plane.}}
\end{figure}

\begin{remarks}{\rm 
  The following facts are easy to see:
\begin{itemize}
  \item [(i)]the domains $D'_\mu$ and $D_\mu$ are biholomorphically equivalent via the mapping 
\begin{align*}
  \varphi:\, &D'_\mu\to
               D_\mu\qquad\qquad\varphi(z_1,z_2):=(e^{z_1}, z_2) \,;
\end{align*}
\item[(ii)] for every $z_1$ fixed the fiber over $z_1$
  $$
  \Pi(z_1)=\{z_2\in\bbC: (z_1,z_2)\in D'_\mu\}
$$
is connected, while  the same property does not hold for $D_\mu$.
\end{itemize}
}
\end{remarks}

For a given domain $\Omega$ in $\bbC^2$ that is rotationally
invariant in
the second variable $z_2$, such as $\cW, D_\mu, D'_\mu$,  using
Fourier expansion in $z_2$ we can decompose the Bergman space
$A^2(\Omega)$ as 
\begin{equation}\label{decom}
A^2(\Omega) =\bigoplus_{j\in\bZ} \cH^j\, , 
\end{equation}
where 
$$
\cH^j =\bigl\{ F\in A^2:\, 
F(z_1,e^{i\theta} z_2) =e^{ij\theta} F(z_1,z_2)\bigr\} \,.
$$

 If for every $z_1$ fixed, the fibers
 $\Pi (z_1)=\{z_2: (z_1,z_2)\in \Omega \}$ are connected, then
$F\in \cH^j $ has the form
$$
F(z_1,z_2) =f(z_1) z_2^j \,,
$$
where $f$ is holomorphic.  In the case of $D'_\mu$, the fibers
$\Pi(z_1)$
are connected and $f$ is holomorphic on the
strip $\{ |\Im z_1|< \mu+\pi/2\}$.  Hence, we may write the kernel $K'$ of $D'_\mu$ as 
$$
K'(\z,\omega)=\sum_{j=-\infty}^{\infty}K'_j(\z_1,\omega_1)\z_2^j\overline\omega_2^j
$$
and using these observation an explicit computation in $1$-dimension, it is possible to compute
the Bergman kernel $K'$ of $D'_\mu$ quite explicitly. In \cite{Ba-Acta} the kernel $K'_{-1}$ is explicitly computed and it holds
$$
K'_{-1}(\z,\omega) = \frac{1}{2\pi} \int_\bbR
\frac{t^2}{\sinh (2\mu t)  \sinh (2\pi t)} e^{i(\z_1-\ov\omega_1)} \, dt \,.
$$
The analysis of the kernels $K'_j$'s for $j\neq-1$ is performed in \cite{KP-Houston} and it is more difficult
since some cancellations that simplify the computations in the case $j=-1$ do not occur for $j\neq -1$.
Recalling the transformation rule for the Bergman kernel,
$$
 P_{D'_\mu}[\varphi' (f\circ\varphi)]= \varphi'
 [(P_{D_\mu}f)\circ\varphi] \,,
$$
Barrett analyses the kernel $K$ of $D_\mu$ and, in particular, the $(-1)$-component of the kernel, concluding that  $P_{D_\mu}$ is not exactly regular, that is, $P_{D_\mu}$ is not a bounded operator $P_{D_\mu}: W^{s,2}(D_\mu)\to W^{s,2}(D_\mu)$ for $s$ sufficiently large. The same conclusion for the smooth worm $\mathcal W$ is then obtained via an exhaustion argument. Setting $\cW^\tau =\big\{(z_1,z_2)\in\bbC^2: (\frac{z_1}{\tau},z_2)\in\cW\big\}$,
Barrett showed that
\begin{equation}\label{projection-limit}
P_{\cW^\tau} f\to P_{D_\mu}f
\end{equation}
as $\tau\to\infty$. For, if we denote by $d_\tau$ the dilation in the
first variable by $\tau>0$, $d_\tau(z_1,z_2)= (\tau z_1,z_2)$, then
$d_\tau(\cW) = \cW^\tau$ and 
\begin{equation}\label{projection-transform}
P_{\cW^\tau}=T^{-1}_\tau P_{\mathcal W} T_\tau.
\end{equation}
where $T_{\tau} f(z_1,z_2)=f(\tau z_1,z_2)$. From this relation, it is possible to deduce the boundedness of $P_{\cW^\tau}$ from the one of $P_\mathcal W$. Then, passing to the limit as in \eqref{projection-limit}, we would obtain the boundedness of $P_{D_\mu}$, hence, a contradiction. In order to prove \eqref{projection-limit}, it is necessary the trivial, but important, remark that given any compact set $E\subseteq
D_\mu$, there exists $\tau_E>0$ so
that for all $\tau\ge \tau_E$, $E\subseteq
\cW^\tau$.

Thus, the analysis on
the domains $D'_\mu$ and $D_\mu$ not only provided intuition on the
case of the smooth, bounded worm domain $\cW$, but also it was
fundamental in proving the result on the irregularity of the Bergman
projection on $\cW$ itself. 
\ms

We now briefly comment on Christ's result \cite{Christ}.  He proved
that $\cW$ does not satisfy Condition (R), by showing that for all
$s>0$
(apart from a discrete set of
exceptions) the Neumann operator $\cN$ satisfies an {\em a priori} estimate 
$$
||\cN u||_{W^{2,s}} \leq C_{s,j} ||u||_{W^{2,s}}
$$
valid for every $u \in \cH^j_1 \cap C^\infty(\overline{\cW})$ 
such that $\cN u \in C^\infty(\overline{\cW})$.  (Here the subscript
$1$ indicates the fact that $u$ is a $(0,1)$-form.)
If $\cN: C^\infty(\overline{\cW})\to C^\infty(\overline{\cW})$ were
bounded,
such estimates 
would 
contradict the irregularity of $P_{\cW}$.

We conclude this section by discussing the Diederich--Forn\ae ss index
of the worm domain $\cW$.  In \cite{Liu1,Liu2} B. Liu proved that
$\operatorname{DF}(\cW)= \pi/2\mu$, see also \cite{K-Liu-P}.

\ms

\section{Hardy spaces on model worm domains}\label{Har-sec}

We now consider another canonical kernel and projection of a domain
$\Omega$ in $\bbC^n$, the Szeg\H o kernel and projection.

Let 
$$
\Omega =\big\{z\in\bbC^n:\, \rho(z)<0 \big\} \,,
$$
where $\rho$ is smooth and $\nabla\rho\neq 0$ on $\p \Omega$.   Let
$\Omega_\eps=\{\rho(z)<-\eps\}$, and suppose there exists a
family of Borel measures $\{\sigma_\eps\}$ on $\ov\Omega$ and
supported 
 on $\p \Omega_\eps\subseteq \ov \Omega$ such 
$\sigma_\eps \to \sigma_0=:\sigma$ weakly as $\eps\to0$, that is, for
all $f\in C(\ov\Omega)$, $\int f d\sigma_\eps \to \int fd\sigma$ as $\eps\to0$.  
Define the Hardy space $H^2(\Omega,d\sigma)$ as  
$$
H^2(\Omega,d\sigma)=\Big\{f\in \operatorname{Hol}(\Omega):
\,
\sup_{\eps>0}\int_{\p\Omega_\eps}|f(\zeta)|^2d\sigma_\eps(\zeta)<\infty\Big\} \,.
$$

Under mild conditions on the family of measures $\{\sigma_\eps\}$, the
Hardy
$H^2(\Omega,d\sigma)$ is a reproducing kernel Hilbert space and
its
reproducing kernel is called the Szeg\H o kernel.  The classical case
is $\Omega$ is a smoothly bounded domain and 
$d\sigma_\eps$ is the induced surface measure on $\p\Omega_\eps$.
In this case, we
simply write $H^2(\Omega)$.  
It is a classical result (see \cite{St2}) that under these assumptions, if
$f\in H^2(\Omega)$ then $f$ converges {\em non-tangentially} to a
boundary function $\widetilde f\in L^2(\p\Omega)$.  In full
generality, these latter facts have to be shown to hold true.

Thus, we may define
$$
H^2(\p\Omega,d\sigma) =\Big\{g\in L^2(\p\Omega,d\sigma): \,
g(\z)=\lim_{z\to\zeta}f(z) \ \text{non-tangentially, for some}\ f\in
H^2(\Omega,d\sigma)\Big\} \,.
$$
Then, the Szeg\H o projection is the Hilbert space orthogonal
projection of $L^2(\p\Omega,d\sigma)$ onto its (closed) subspace $H^2(\p\Omega,d\sigma)
$, the subspace of boundary
values of functions in $H^2(\Omega,d\sigma)$,
$$
S_\Omega : L^2(\p\Omega,d\sigma)\to H^2(\p\Omega,d\sigma)\qquad
S_\Omega g(\zeta)=\lim_{z\to\zeta\in
  \p\Omega}\int_{b\Omega}g(\zeta')K(z,\zeta') \,d\sigma(\zeta') \,.
$$

By definition, the  Szeg\H{o} projection depends
on the choice of the measure on the boundary.   Another very natural,
and thus far little considered, possible choice, is the
Fefferman surface measure $\sigma_F$, see
 \cite{MR3145917}, or any other surface measure  $\omega
 d\sigma$, where  $\omega$ is a continuous positive function on
 $\p\Omega$, 
 \cite{2015arXiv150603965L}.  The surface measure $\sigma_F$ was introduced
by C. Fefferman in order to obtain a measure that is {\em
  biholomorphic invariant}.  To be precise, suppose $\Omega_1$ and
$\Omega_2$ are bounded, smooth, pseudoconvex domains that admit 
a biholomorphic map $\vp:\Omega_1\to\Omega_2$ that extends to a smooth
$C^\infty$-diffeormophism of the boundary, such as in case one of the
two domains satisfies Condition (R).  Then, the mapping
$\Lambda(f):=\sqrt{\operatorname{det}\vp'} (f\circ\vp)$ defines an isometric
isomorphism $\Lambda: H^2(\Omega_2, d\sigma_F) \to H^2(\Omega_1,
d\sigma_F)$.

 When we consider non-smooth domains, such as $D_\mu$ and $D'_\mu$,
 and more noticeably
 the
 polydisk, it is  perhaps more natural, and certainly interesting to study
 Hardy spaces defined by integration over the so-called {\em
   distinguished boundary}.  Given a domain
 $\Omega\subseteq\bbC^n$, we call the
 distinguished boundary, and we denote it by $d_b(\Omega)$,the set
\begin{equation}\label{dist-bndry-def1}
d_b(\Omega)=\Big\{\zeta\in \p\Omega:  \sup_{z\in \p\Omega} |f(z)| \le
\sup_{\zeta\in d_b(\Omega)} |f(\zeta)| \ \text{for all\ } f\in H^\infty(\Omega)\Big\} 
\end{equation}
  where $H^\infty(\Omega)$ is the spaces of  holomorphic
  functions on $\Omega$ that are bounded.  Then, we consider the
  induced measure on  $d_b(\Omega)$ and denote it by $d\beta$. \ms

  We now describe the main results that we have obtained on the
  regularity of Szeg\H o projections on model worm domains.  We first
  consider the case of $D'_\mu$ and the induced surface measure
  $d\beta$. The following result is in \cite{MP}.

  \begin{thm}\label{Sze-pro-main-thm}
The Szeg\H o projection $\cS$, initially defined
on the dense subspace 
$W^{s,p}(\p D'_\mu)\cap L^2( \p D'_\mu,d\sigma)$, 
 extends to a bounded operator 
$$
\cS: W^{s,p}(\p D'_\mu)\to
 W^{s,p}(\p D'_\mu)\,, 
$$
for $1<p<\infty$ and $s\ge0$.
\end{thm}

The proof of such result relies on explicit computations on the
boundary of $D'_\mu$, which can be written as union of four pieces
that have intersection of null measure. The 
Szeg\H o projection can be correspondingly written as sum of 16
different integral operators.  For each of these operators we apply a
decomposition similar to \eqref{decom} and obtain an explicit
expression and thus write them as composition of a bounded Fourier
multiplier and an operator of Hilbert-type.

It is worth to remark that the boundedness of the corresponding
Szeg\H o projection on $D_\mu$ is still unexplored and it would be
significant to study such (ir)-regularity.
\ms

We now turn to the case of the Szeg\H o projection on the distinguished boundaries (\cite{monguzzi-concrete}). More precisely, denote by $d_b (D'_\mu) $
and $d_b (D_\mu)$ the distinguished boundaries of the domains $D'_\mu $
and $ D_\mu$, resp., and by $\sS'$ and $\sS$ the corresponding
Szeg\H o
projections, resp.  We point out that in this setting, the operators
$\sS'$ and $\sS$ are given by singular
integrals over  $d_b(D'_\mu)$
and $d_b(D_\mu)$, resp.

The
case of $D'_\mu$ was considered in \cite{M}, where the main result is
the following
 \begin{thm}\label{Sze-pro-dist-bndy-thm-D-prime}
The Szeg\H o projection $\sS'$, initially defined
on the dense subspace 
$W^{s,p}(d_b(D'_\mu),d\beta)\cap L^2( d_b(D'_\mu),d\beta)$, 
 extends to a bounded operator 
$$
\sS': W^{s,p}(d_b(D'_\mu),d\beta)\to
 W^{s,p}(d_b (D'_\mu),d\beta)\,, 
$$
for $1<p<\infty$ and $s\ge0$.
\end{thm}
\ms

The proof of this result follows from explicit computations of the Szeg\H o projection of suitably defined Hardy spaces on the distinguished boundary of $D'_\mu$. Such a boundary is the union of four different connected components which are mutually disjoint and the Szeg\H o projection turns out to be a linear combination of bounded Fourier multiplier operators. A detailed analysis of the Szeg\H o kernel associated to $\sS'$ is performed in \cite{monguzzi-CAOT}.

With a similar proof the analogous result for the Szeg\H o projection $\sS$ on the
distinguished boundary of $D_\mu$ is studied and we now describe it with greater details. 
For $(t,s)\in (0,\pt)\times[0,\mu)$ consider the 
domain
$$
D_{t,s}=\left\{(z_1,z_2)\in\bbC^2: \big|\arg\, z_1-\log|z_2|^2\big|<t,
  \big|\log|z_2|^2\big|<s\right\}. 
$$
Then, the domains $\{D_{t,s}\}_{t,s}$ constitute a family of
approximating domains for $D_\mu$. The distinguished boundary
of these domains is given by
\begin{equation*}
d_b(D_{t,s})=\left\{(z_1,z_2)\in\bbC^2: \big|\arg\,
  z_1-\log|z_2|^2\big|=t, \big|\log|z_2|^2\big|=s\right\}. 
\end{equation*}
Consequently,
for $1\le p<\infty$, we define the Hardy space $H^p(d_b(D_\beta),d\beta)$
defined by
\begin{align*}
  H^p(d_b(D_\beta),d\beta)
  &  =\bigg\{ f\in \Hol(D_\beta):\, \|f\|^p_{H^p(d_b(D_\beta),d\beta)}
    = \sup_{(t,s)\in(0,\pt)\times[0,\mu)}
\|f\|^p_{L^p(d_b(D_{t,s}),d\beta)}<\infty\bigg\} \,,
\end{align*}
where, denoting by  $d\beta_{t,s} $ the induced measure on $d_b(D_{t,s})$,
\begin{align*}
&  \|f\|^p_{L^p(d_b(D_{t,s}),d\beta)}
= \int_{d_b(D_{t,s})} |f|^p\,
 d\beta_{t,s} \notag \\
& =  \int_0^{\infty}\int_0^{2\pi}|f\big(re^{i(s+t)}, e^{\frac{s}{2}}
e^{i\theta}\big)|^p\ e^{\frac{s}{2}} d\theta dr  +\int_0^{\infty}\int_0^{2\pi} |f\big(r e^{i(s-t)},
 e^{\frac{s}{2}}e^{ i\theta}\big)|^p \ e^{\frac{s}{2}} d\theta dr
  \\ 
 & \qquad +\int_0^{\infty}\int_0^{2\pi} |f\big(re^{-i(s+t)},
 e^{-\frac{s}{2}}
e^{i\theta}\big)|^p\ e^{-\frac{s}{2}} d\theta dr
+\int_0^{\infty}\int_0^{2\pi} |f\big(r e^{-i(s-t)},
  e^{-\frac{s}{2}}e^{ i\theta}\big)|^p \ e^{-\frac{s}{2}}d\theta dr \,.\notag
\end{align*}

The main results in \cite{MP2} 
are the following.
The first result provides the sharp interval of values of $p$ for which
the Szeg\H o projection $\sS$ on the distinguished boundary of $D_\mu$
is bounded.  We recall that we set $\nu=\frac{\pi}{2\mu}$, so that
$\nu=\nu_\mu$ tends to $0$ as $\mu$ becomes large.

\begin{thm}\label{t:LpBounds}
The Szeg\H{o}
projection $\sS$, initially defined on the dense subspace
$L^p(d_b(D_\mu),d\beta) \cap L^2(d_b(D_\mu),d\beta))$,  extends to a bounded operator   
$$
\sS:L^p(d_b(D_\mu),d\beta))\to L^p(d_b(D_\mu),d\beta))
$$
if and only if $\frac{2}{1+\nu}<p<\frac{2}{1-\nu}$.
\end{thm}

The next result concerns with the sharp boundedness of $\sS$ on the
$L^2$-Sobolev spaces on $d_b(D_\mu)$. 

\begin{thm}\label{SobolevL2}
The Szeg\H{o}
projection $\sS$ defines a bounded operator  
$$
\sS: W^{s,2}(d_b(D_\mu),d\beta))\to W^{s,2}(d_b(D_\mu),d\beta)
$$
if and only if $0\leq s<\frac{\nu}{2}$.
\end{thm}

In the case of Sobolev norms with $p\neq 2$ we do not have a complete
characterization of the mapping properties of $\sS$, but we have a
partial result.
\begin{thm}\label{SobolevLp}
 Let $s>0$ and
 $p\in(1,\infty)$. If the operator $\sS$, initially defined on the
 dense subspace $W^{s,p}(d_b(D_\mu),d\beta))\cap L^2(d_b(D_\mu),d\beta)$, extends to
 a bounded operator $\sS: W^{s,p}(d_b(D_\mu),d\beta)\to W^{s,p}(d_b(D_\mu),d\beta)$, then 
$$
-\frac{\nu_\beta}{2}\leq s+\frac{1}{2}-\frac{1}{p}\leq\frac{\nu_\beta}{2}.
$$
Assuming $p\geq2$ we obtain the stronger condition
\begin{equation*}
 0\leq s+\frac12-\frac1p<\frac{\nu_\beta}{2}.
\end{equation*}
\end{thm}

The main fact used in the proofs is that we can write
the Szeg\H{o} projection $\sS$ as a sum of
Mellin--Fourier multiplier operators which we now briefly describe. In order to do so we introduce
some notation. We set $\bbX$ to denote either
$\bbR$ or $\bbT$,  and, accordingly, 
$\widehat{\bbX}= \bbR$, or $\bbZ$, respectively, where we denote by
$\widehat{\,}$ the Fourier transform or Fourier series on $\bbR$ and
$\bbT$, resp.
Instead, we denote by $\cF$ the Fourier transform on $\bbR\times\bbX$, given by
$$
\cF f(\xi_1,\xi_2) =
\int_{\bbR\times\bbX} f(x_1,x_2) e^{- i(x_1\xi_+x_2\xi_2)}\, dx_1dx_2
$$
when $f$ is absolutely integrable.
We consider
   the Fourier
multiplier operator given by
$$
T_m(f) = \cF^{-1} \big(m \cF f\,\big) \,
$$
when $m$ is a bounded measurable function on
$\bbR\times\widehat{\bbX}$.
We say that a bounded function $m$ on $\bbR\times\bbX$ is a {\em
  bounded Fourier multiplier} on $L^p(\bbR\times\bbX)$ if 
$T_m: L^p(\bbR\times\bbX)\to L^p(\bbR\times\bbX)$ is bounded.

Given a function  $\vp\in C^\infty_c \big((0,\infty)\times
\bbX\big)$ we define the operator 
 $$
\cC_{p} \vp(x,y) = e^{\frac 1p(x)} \vp(e^{x},y) \,.
$$
It is clear that $\cC_{p}$ extends to an isometry of  $L^p\big((0,\infty)\times
\bbX\big)$  onto
$L^p \big((0,\infty) \times
\bbX\big)$.  

For $a,b\in\bbR$, with $0<a<b<1$, we denote by $S_{a,b}$ the vertical
strip in the complex plane
\begin{equation*}
S_{a,b}=\big\{z\in\bbC:\,  a<\Re z<b \big\}\,.
\end{equation*}
Given a bounded measurable function $m$ defined on
$S_{a,b}\times\bbX$, when $a<\frac1p<b$
we write
\begin{equation*}
m_p(\xi_1,\xi_2)=m(\textstyle{\frac1p}- i\xi_1,\xi_2) \,. 
\end{equation*}

Finally,
 we define an operator acting on functions defined on 
$(0,+\infty)\times \bbX$ 
as
\begin{equation}\label{def-MF-multiplier}
\sT_{m,p} =\cC_p^{-1} T_{m_p} \cC_p \,.
\end{equation}
We call such an operator a  Mellin--Fourier multiplier operator, the
reason for which will soon be  clear.  A similar
class of operators was studied by Rooney \cite{MR860095}.  
Incidentally, we believe that this class of operators is of its own interest. 

\begin{thm}\label{rooney-thm}
With the above notation, let $m:S_{a,b}\times\bbX\to \bbC$ be
continuous and such that 
\begin{itemize}
\item[(i)]
$m(\cdot,\xi_2)\in\Hol(S_{a,b})$ and bounded in every closed substrip
of $S_{a,b}$,
for every $\xi_2\in\bbX$ fixed;\smallskip
\item[(ii)] for every $q$ such that $a<\frac1q<b$,
$m_q$ is a
  bounded Fourier multiplier on $L^q(\bbR\times\bbX)$.
\end{itemize}
Then, for $a<\frac1p<b$,  $\sT_{m,p}=\sT_m$ is independent  
of $p$ and  
$$
\sT_m: L^p((0,+\infty)\times\bbX)
\to L^p( (0,+\infty)\times\bbX) 
$$
is bounded.
\end{thm}

In the course of the proof we show that if $m$ satisfies the hypotheses of
the theorem, then
\begin{equation}\label{M-F-mult-eq}
\cC_p^{-1}  T_{m_p}\cC_p (\vp)
= \cF_2^{-1} M_1^{-1} \big( 
m
( M_1 \cF_2\vp) \big) \,,
\end{equation}
where $M_1$ denotes the Mellin transform in the first variable, that is,
$$
M_1\varphi (z,\z_2)=\int_0^{+\infty}t^{z-1}\varphi(t,\zeta_2)\, dt,
$$
and $\mathcal F_2$ denotes the Fourier transform in the second variable.
Equality \eqref{M-F-mult-eq} clearly giustify the fact that 
the operator $\sT_m$ a Mellin--Fourier multiplier operator:
 it is a Mellin transformation in the first variable, a Fourier
 transformation in the second variable, followed by multiplication by
 $m$ and then the inverses of the Mellin and Fourier transforms.
We  also point that, if $m,\tilde m$ satisfy the assumptions in
the theorem, then $\sT_m \sT_{\tilde m}= \sT_{m\tilde m}$, and thus it
is reasonable to call these operators {\em multipliers}. 

Once we explictly write the Szeg\H o projection $\sS$ as a linear
combination of Mellin--Fourier multiplier operators, we are able to
study its regularity by also exploiting the regularity of $\sS'$. We
recall that, unlike in the case of the Bergman projection, 
in the Szeg\H{o} setting, in general, there is no transformation rule
for the Szeg\H{o} projection under biholomorphic mappings. 
Nonetheless, we are able to prove a transformation rule for the
projections $\sS$ and $\sS'$.  Recall that $D'_\mu$ and $D_\mu$ are biholomorphically equivalent via the map
\begin{align}\label{Biholomorphism}
\begin{split}
\vp^{-1}: &D_\beta\to D'_\beta\\
 &(z_1,z_2)\mapsto (\Log(z_1 e^{-i\log|z_2|^2})+i\log|z_2|^2, z_2),
\end{split}
\end{align}

where $\Log$ denotes the principal branch of the complex logarithm. Setting 
\begin{equation*}
 \psi_p(z_1,z_2) := e^{-\frac ip\log|z_2|^2}(z_1
 e^{-i\log|z_2|^2})^{-\frac 1p} 
\end{equation*}
we obtain that 
$$
\sS'(\Lambda^{-1}f)=\Lambda^{-1}(\sS),
$$
where $\Lambda f:=\psi_p(f\circ\varphi^{-1})$.
\ms

%
%
%
%
%
%

\section{Other results on the regularity of Szeg\H o projections.}
\label{Other-section}

Mapping properties of the Szeg\H o projection on other function spaces
have been studied for various classes of smooth bounded
domains and the are several positive results. The Szeg\H o projection $S_\Omega$ turns out to be bounded on the Lebesgue--Sobolev spaces
$W^{s,p}(\p\Omega)$ for $1<p<\infty$ and
$s\ge0$ in the case of strictly pseudoconvex domains
\cite{MR0450623},  domains of finite type in $\bbC^2$
\cite{MR979602} and convex domains of finite type in
$\bbC^n$ \cite{MR1452048}. The exact regularity of $S_\Omega$, that is, the boundedness $S_\Omega: W^{s,2}(\p\Omega)\to W^{s,2}(\p\Omega)$
for every $s\geq0$, holds when $\Omega$ is a Reinhardt domain \cite{MR773403,MR835396},
a domain with partially transverse symmetries \cite{MR971689}, a
pseudoconvex domain satisfying Catlin's property $(\cP)$ \cite{
  MR871667}, a 
complete Hartogs domain in $\bbC^2$ \cite{ MR999739}, or a domain with
a plurisubharmonic defining function on the boundary \cite{MR1133741}.   We also mention that, if $\Omega$ is bounded,
$C^2$ and strongly pseudoconvex in $\bbC^n$, the Szeg\H{o} projection
$P_\Omega$ again extends to bounded operator on $L^p(\p\Omega)$ for
$1<p<\infty$, 
\cite{2015arXiv150603748L,2015arXiv150603965L}. 

There are also examples of domains $\Omega$  on which the Szeg\H{o} projection
$P_\Omega$ is less regular.  L. Lanzani and E. M. Stein described the
(ir-)regularity of $P_\Omega$ on Lebesgue spaces in the case of planar
simply connected domains, \cite[Thm. 2.1]{MR2030575}.  In particular
they showed that if $\Omega$ has Lipschitz boundary, then
$P_\Omega:L^p(\p\Omega)\to L^p(\p\Omega)$ if and only if
$p'_\Omega<p<p_\Omega$, where $p_\Omega$ depends only on the Lipschitz
constant of $\p\Omega$. 
More recently, S. Munasinghe and Y.E. Zeytuncu provided an example of
a piecewise smooth, bounded pseudoconvex domain in $\bbC^2$ on which
the Szeg\H{o} projection
$P_\Omega$ is unbounded on
$L^p(\p\Omega)$ for every $p\neq 2$ \cite{MR3355787}. The same result on tube
domains over irreducible self-dual cones of rank greater than 1 has
been known for a number of years, \cite{BeBo}. 

In a recent paper \cite{LS-worm} Lanzani and Stein announced a result concerning the $L^p$ continuity of the Szeg\H o projection  attached to the smooth worm domain $\mathcal W_\mu$ with respect to the induced surface measure d$\sigma$ on $\partial \mathcal W_\mu$.
In particular, they announced that for any $p\neq 2$ there is a $\mu=\mu(p)$ such that the Szeg\H o projection  is not bounded $P_{\mathcal W_\mu}:L^p(\partial\mathcal W_\mu)\to L^p(\partial\mathcal W_\mu)$. 

It is reasonable to think that the culprit of the (ir-)regularity of both the Bergman and Szeg\H o projection on the worm domain $\mathcal W_\mu$ is the presence of the critical annulus $\cA= \{ (0,z_2):\, \big|\log |z_2|^2\big|\le
\mu\}$ in the boundary $\partial \mathcal W_\mu$; see, for instance,
\cite{BS-worm}. For this reason, it would be interesting to study the
Szeg\H o projection of $\mathcal W_\mu$ with respect to the Fefferman
measure $d\sigma_F$. In fact, as we now see,
the Fefferman
measure $d\sigma_F$ is given by a smooth density $\omega$ times
$d\sigma$ and the density $\omega$ vanishes identically on
the critical annulus $\cA$/ In detail, given
$\Omega=\{z\in\bbC^n:\rho(z)<0\}$, the Fefferman surface area measure 
(\cite[pg. 259]{Fefferman},\cite{MR3145917})
on $\partial\Omega$ is defined by 
\begin{equation}\label{fefferman-measure}
 d\sigma_F=c_n\sqrt[n+1]{M(\rho)}
 \frac{d\sigma}{\|\nabla \rho\|}
\end{equation}
where $M(\rho)$ is the Fefferman Monge--Amp\'ere operator 
$$
M(\rho)={-\text{det}
 \begin{pmatrix}
\rho && \rho_{\overline k}\\
\rho_j && \rho_{j\overline{k}}
\end{pmatrix}_{1\leq j,k\leq n}
}.
$$
It can be proved (see also \cite[Section 2]{MR3145917}) that the definition of $d\sigma_F$ does not depend on the defining function $\rho$ and that there exists a sesqui-holomorphic kernel $S(z,\zeta)$ such that, for every $f\in H^2(\Omega)$,
$$
f(z)=\int_{\partial\Omega} f(\zeta) S(z,\zeta)\, d\sigma_F(\zeta).
$$
Hence, Hardy spaces and Szeg\H o projections with respect to the Fefferman measure can be defined and investigated. 
%
In the case of the worm domain $\mathcal W_\mu$ the defining function is
\begin{align*}
\rho(z_1,z_2)&= |z_1|^2-2\Re(z_1 e^{-i\log|z_2|^2})+\eta(\log|z_2|^2),
\end{align*}
therefore, setting $R(z_1,z_2) = \Re(iz_1 e^{-i\log|z_2|^2})$,. 
\begin{align*}
M&(\rho)(z_1,z_2)=\\
&\begin{pmatrix}
0 && z_1-e^{i\log|z_2|^2} && \frac{1}{\overline{z_2}}\big(2R(z_1,z_2)+\eta'(\log|z_2|^2)\big)\\ \vspace{2pt}
\overline{z_1}-e^{-i\log|z_2|^2} && 1 && \frac{i}{\overline{z_2}} e^{-i\log|z_2|^2}\\ \vspace{2pt}
\frac{1}{z_2}\big(2R(z_1,z_2)+\eta'(\log|z_2|^2)\big) && -\frac{i}{z_2} e^{i\log|z_2|^2} && \frac{1}{|z_2|^2}\big(2R(z_1,z_2)+\eta''(\log|z_2|^2)\big)
\end{pmatrix}.
\end{align*}
When we restrict the matrix $M(\rho)$ to the critical annulus $\cA$ we get 
$$
\det M(\rho)(0,z_2)=\det\begin{pmatrix}
 0 && -e^{i\log|z_2|^2} && 0\\
 -e^{-i\log|z_2|^2} && 1 && \frac{i}{\overline{z_2}}e^{-i\log|z_2|^2}\\
 0 && -\frac{i}{z_2}e^{i\log|z_2|^2} && 0
\end{pmatrix}=0.
$$

Since the boundary of the domain $D_\mu$ (similarly, of $D'_\mu)$ is Levi flat, that is, the Levi form of its defining function is identically zero at every point of $bD_\mu$, the density of the Fefferman measure on $bD_\mu$ is identically zero. This can be easily verified by explicitly computing $M(\rho)$ for $\rho(z_1,z_2)=\Re(z_1 e^{-i\log|z_2|^2}).
$
Thus, the Szeg\H o projection on $\cW_\mu$ with respect to the Fefferman area measure cannot be investigated exploiting the model domains, but it must be directly approached. This certainly is an interesting direction for future research.
\ms

\section{Hartogs triangles}\label{Hartogs-sec}

Another class of domains on which it is interesting to test and study the regularity of the Bergman and Szeg\H o projection is the one of generalized Hartogs triangles.  
Given a real parameter $\gamma>0$, the generalized Hartogs triangle $\mathbb H_\gamma$ is defined as 
$$
\mathbb H_\gamma=\big\{(z_1,z_2)\in\bbC^2: |z_1|^\gamma<|z_2|<1\big\}.
$$
This family of domains were recently introduced in \cite{E-pacific} and the value $\gamma=1$ corresponds to the classical Hartogs triangle $\mathbb H$ (\cite{Shaw}). The Hartogs triangle is a simple, but not trivial, model domain on which it is worth to test several conjectures. It turns out that $\mathbb H$ is a source of counterexamples in complex analysis. For instance, as the worm domain $\mathcal W$, the Hartogs triangle has non-trivial Nebenh\"ulle. However, unlike $\mathcal W$, the domain $\mathbb H$ is not smooth; on the contrary, it is highly singular at the point $z_1=z_2=0$.  This pathological geometry affects the $L^p$ behavior of the Bergman projection $P_{\mathbb H}$ and it turns out that $P_\mathbb H$ extends to a bounded operator $P_\mathbb H:L^p(\mathbb H)\to L^p(\mathbb H)$ if and only if $p\in(\frac43,4)$ (\cite{CZ}). This result has been extended to the case of generalized Hartogs triangle in a series of paper by L. D. Edholm and J. D. McNeal and it holds that $P_{\mathbb H_\gamma}$ extends to a bounded operator $L^p(\mathbb H_\gamma)\to L^p(\mathbb H_\gamma)$ for a restricted range of $p\in(1,\infty)$ whenever $\gamma\in\mathbb Q$, but $P_{\mathbb H_\gamma}$ is unbounded on $L^p(\mathbb H_\gamma)$ for any $p\neq2$ whenever $\gamma$ is irrational (\cite{E-pacific, EMcN, EMcN17}). The Sobolev (ir-)regulariy of $P_{\mathbb H_\gamma}$ has been investigated as well and we refer the reader to the very recent paper \cite{EMcN-Sobolev}.

In addition to the aforementioned papers, we mention also the recent papers \cite{HW, HW2}, where weighted $L^p$ and endpoint estimates for the classical Hartogs triangle are obtained via dyadic harmonic analysis techniques, and  \cite{Chen5,Chen1,Chen2,Chen3,Chen4}, where some other generalizations of the Hartogs triangle and the associated weighted Bergman projections are investigated. 

The definition of a Hardy space $H^2$ on $\mathbb H$, hence the definition of a Szeg\H o projection on $\mathbb H$, is not canonical due to the geometry of the domain. We now recall the definition of a candidate Hardy space on the classical Hartogs triangle which is introduced by the first author in a recent paper \cite{monguzzi-hartogs}. 

Let $\nu>-1$ be a real parameter, let $\bbD$ be the unit disc in the complex plane and let us consider the classical weighted Bergman spaces $A^2_\nu(\bbD$) defined as the space of holomorphic functions in $\bbD$ endowed with the norm 
$$
\|f\|^2_{A^2_\nu(\bbD)}=(\nu+1)\int_{\bbD}|f(z)|^2(1-|z^2|)^\nu\, dz.
$$
It is a well-known fact that 
$$
\|f\|^2_{H^2(\bbD)}=\lim_{\nu\to-1^+}\|f\|^2_{A^2_\nu(\bbD)}
$$
where $H^2(\bbD)$ is the Hardy space in $\bbD$, that is, the space of holomorphic functions in $\bbD$ endowed with the norm
$$
\|f\|^2_{H^2(\bbD)}:=\sup_{0<r<1}\frac{1}{2\pi}\int_0^{2\pi}|f(re^{i\theta})|^2\, d\theta.
$$
Notice that if $K_{\bbD}(z,w)=(1-\overline z w)^{-2}$ denotes the reproducing kernel of $A^2(\mathbb H)$, the unweighted Bergman space, then 
$$
K^{-\frac\nu2}_\bbD(z,z)=(1-|z|^2)^{\nu}\approx\delta^\nu(z),
$$
where $\delta(z)$ is the distance of $z\in\bbD$ from the topological boundary $\partial \bbD$. Therefore, we analogously define the weighted Bergman space $A^2_\nu(\mathbb H)$ as the space of holomorphic functions on $\mathbb H$ such that 
\begin{equation}\label{norm-bergam-hartogs}
\|f\|^2_{A^2_\nu(\mathbb H)}:= C_\nu \int_\mathbb H |f(z)|^2 K^{-\frac\nu2}(z,z)\, dz,
\end{equation}
where $C_\nu$ is a positive constant to be chosen, $dz$ denotes the Lebesgue measure in $\mathbb C^2$ and $K(z,w)$ is the reproducing kernel of the unweighted Bergman space $A^2(\mathbb H)$ and it is given by
$$
K(z,w)=K\big((z_1,z_2),(w_1,w_2)\big)=\frac{1}{2z_2\overline w_2(1-\frac{z_1\overline w_1}{z_2\overline w_2})(1-z_2\overline w_2)^2}.
$$

The following proposition is proved in \cite[Theorems 3.1.4 and 3.1.5]{E-thesis} and describes the diagonal behavior of the kernel $K$.
\begin{prop}[\cite{E-thesis}] The following facts hold true.
\begin{enumerate}
             \item[$(i)$] Let $\delta(z)$ be the distance  of $z$ to $\partial\bbH$, the topological boundary of $\bbH$. Then, 
             $$K(z,z)\approx \delta(z)^{-2}
             $$
             as $z$ tends to the origin.
\medskip
\item[$(ii)$] Let $p$ be any point in the distinguished boundary $d_b(\bbH)$.   For any number $\beta\in (2,4]$ there exists a path $\gamma: [1/2,1]\to\overline\bbH$ such that $\gamma(1)=p$ and for all $u\in[1/2,1)$,
$$
K(\gamma(u),\gamma(u))\approx \delta(\gamma(u))^{-\beta}.
$$
\end{enumerate}
\end{prop}
In \cite{monguzzi-hartogs} the $L^p$ regularity of the weighted Bergman projection $P_\nu$ is completely characterized. 
\begin{thm}[\cite{monguzzi-hartogs}, Theorem 1]\label{monguzzi-1}
Let $\nu>-1$ and let $P_\nu$ be the weighted Bergman projection densely defined on  $L^2_\nu(\bbH)\cap L^p_\nu(\bbH)$ for $p\in (1,+\infty)$. Then, we have the following:
 \begin{enumerate}
 \item[(i)] if $\nu>0$ and $\nu\neq 2n,n\in\mathbb N$, the weighted Bergman projection $P_\nu$ extends to a bounded operator $P_\nu: L^p_\nu(\bbH)\to L^p_\nu(\bbH)$ if and only if $p\in \big(2-\frac{\nu-2[\nu/2]}{2+\nu-[\nu/2]}, 2+\frac{\nu-2[\nu/2]}{2+[\nu/2]}\big)$;
 \medskip
 
 \item[(ii)] if $\nu=2n,n \in\mathbb N_0$, the weighted Bergman projection $P_\nu$ extends to a bounded operator $P_\nu: L^p_\nu(\bbH)\to L^p_\nu(\bbH)$ if and only if  $p\in \Big(2-\frac{2}{3+n}, 2+\frac 2{1+n}\Big)$;
 
 \medskip
 
 \item[(iii)] if $-1<\nu<0$,  the weighted Bergman projection $P_\nu$ extends to a bounded operator $P_\nu: L^p_\nu(\bbH)\to L^p_\nu(\bbH)$ if and only if $p\in \Big(2-\frac{2+\nu}{3+\nu}, 4+\nu\Big)$.
\end{enumerate}
\end{thm}

The proof of this result follows from an explicit computation of the weighted kernel $K_\nu$ and an application of classical the Schur's lemma to the operator with positive kernel $|K_\nu|$. 

The Hardy space $H^2(\bbH)$ is then defined as the limit space corresponding to the value $\nu=-1$ of the parameter. In particular, $H^2(\bbH)$ is defined in a way such that 
$$
K_{H^2(\bbH)}(z,w)=\lim_{\nu\to-1^{+}}K_{\nu}(z,w)
$$
and
$$
\|f\|^2_{H^2}=\lim_{\nu\to-1^{+}}\|f\|^2_{A^2_\nu}.
$$
It turns out that 
$$
H^2(\bbH):=\bigg\{f\in\Hol(\bbH):\sup_{(s,t)\in(0,1)\times(0,1)}\frac{1}{4\pi^2}\int_{d_b(\bbH_{st})}|f|^2\, d\sigma_{st}<+\infty\bigg\},
$$
where $d\sigma_{st}$ denotes the induced surface measure on $d_b(\bbH_{st})$, the distinguished boundary of the domain 
$$
\bbH_{st}=\left\{(z_1,z_2)\in\bbC^2:|z_1|/s<|z_2|<t\right\}\subsetneq\bbH
$$
for $(s,t)\in(0,1)\times(0,1)$. In particular, $d_b(\bbH_{st})=\big\{(z_1,z_2)\in\bbC^2: |z_1|/s=|z_2|=t\big\}$. We endow $H^2(\bbH)$ with the norm
$$
\|f\|^2_{H^2}:=\sup_{(s,t)\in(0,1)\times(0,1)}\frac{1}{4\pi^2}\int_{d_b(\bbH_{st})}|f|^2\, d\sigma_{st}=\sup_{(s,t)\in(0,1)\times(0,1)} \frac{1}{4\pi^2}\int_0^{2\pi}\int_0^{2\pi}|f(st e^{i\theta}, t e^{i\gamma})|^2 st^2\, d\theta d\gamma.
$$
The Hardy space $H^2(\bbH)$ can be identified with a closed subspace
of $L^2(d_b(\bbH))$, which we denote by $H^2(d_b(\bbH))$, hence a
Szeg\H o projection $S:L^2(d_b(\bbH))\to H^2(d_b(\bbH))$ is
well-defined and can be investigated. In particular the following
holds.  
\begin{thm}[\cite{monguzzi-hartogs}, Theorem 2]
 The Szeg\H o projection $S$ densely defined on $L^2(d_b(\bbH))\cap
 L^p(d_b(\bbH))$ extends to a bounded operator $S:L^p(d_b(\bbH))\to
 L^p(d_b(\bbH))$ for any $p\in(1,+\infty)$. 
\end{thm}
In comparison with Theorem \ref{monguzzi-1}, the $L^p$ regularity of
the Szeg\H o projection is surprising and unexpected. The reason of
this result may be found in the fact that the Hardy space considered,
even if it is naturally defined, turns out to be modeled only on the
distinguished boundary $d_b(\bbH)$ of $\bbH$ and not on the whole
topological boundary $\partial \bbH$. Therefore, we loose track of the
origin $(0,0)$, the most pathological point of $\partial \bbH$.  

A further investigation of Hardy spaces on $\bbH$ and the extension of
the results in \cite{monguzzi-hartogs} to the case of generalized
Hartogs triangle certainly is an interesting direction for future
research. \ms

\section{Orthogonal sets and the {M{\"u}ntz-Sz{\'a}sz} problem for the
  Bergman space}\label{MS-sec}

It would be ideal to be able to obtain the asymptotic expansion of the
Bergman and Szeg\H o kernels on the worm domain $\cW$.  A fundamental
step in this direction would be to obtain the 
explicit expression of Bergman and Szeg\H o kernel on the {\em
  truncated} worm domain
\begin{equation}\label{truncated-worm}
  \cW' = \left\{(z_1,z_2)\in\bbC^2: |z_1-e^{i\log
    |z_2|^2}|^2<1,\ |\log|z_2|^2|<\mu \right\}  \,.
  \end{equation}
Obviously, this domain coincides with $\cW$ when we select
$\eta=\chi_{|t|>\mu}$, and $\cW'$ is bounded, non-smooth, and its
boundary contains the same critical annulus $\cA$ as $\cW$.   In
analogy with the case of the unit bidisk $\bbD^2$, we are led to look for
an orthonormal basis of mononials.  In the case of $\cW'$, as well as
of $\cW$,   the following functions resemble the monomials
$z_1^jz_2^k$, $j,k\in\bbN$, where $\bbN$ denotes the set of
non-negative integers.  We set
$$ 
E_\eta(z) = e^{\eta L(z)}   \, ,
$$
where
$$
L(z) = \log\big(z_1e^{-i\log|z_2|^2}\big)+i\log|z_2|^2\, ,
$$
and $\log$ denotes the principal branch of the logarithm, so that
$$
E_\eta(z_1,z_2) = \big(z_1e^{-i\log|z_2|^2}\big)^\eta
e^{i\eta\log|z_2|^2}\, .
$$
Now we define constants $\gamma_{\alpha\beta}
= h(\alpha-\overline\beta)$, where 
$h(z)=\frac{\sinh
  [\mu(j+1+iz)]}{j+1+iz}$. The following is Proposition 3.1 in \cite{KPS2}.

\begin{prop}\label{prop1}  
 Let $\mu>0$.  
For $\alpha\in\bbC$ and $j\in\bbZ$ let 
$F_{\alpha,j}(z_1,z_2)=E_\alpha(z) z_2^j$.  Then $F_{\alpha,j}\in A^2(\cW'_\mu)$ if
and only if $\Re\alpha>-1$.  Moreover, if $\Re\alpha,\Re\beta>-1$
then
$$
\la F_{\alpha,j},\, F_{\beta,j}\ra_{A^2(\cW'_\mu)}
= (2\pi)^2 \gamma_{\alpha\beta} \frac{\Gamma(\alpha+\overline\beta+2)}{
\Gamma(\alpha+2)\Gamma(\overline\beta+2)}\, .
$$
In particular, 
$\la F_{\alpha,j},\, F_{\beta,j}\ra_{A^2(\cW'_\mu)} =0$
if and only if
\begin{equation}\label{orth-cond}
\alpha-\overline\beta= 2k\nu +i(j+1)\qquad
\text{with\ } k\in\bbZ\setminus\{0\}\, .
\end{equation}
\end{prop}
Thus, if $c>-1$ and $\ell\in\bbN$, and we set
\begin{equation}\label{H-ell-j-def}
H _{\ell,j}(z_1,z_2) =
E_{c_0+\nu\ell+i(j+1)/2}(z)z_2^j\,,
\end{equation}
the next corollary follows.
\begin{cor}\label{ort-sys}  
Each of the two sets
\begin{equation}\label{two-sets}
\big\{  H_{2k,j}\, ,\, j\in\bbZ,\,  k\in\bbN\big\},
\quad\text{and}\quad
\big\{  H_{2k+1,j}\, ,\, j\in\bbZ,\,  k\in\bbN\big\},
\end{equation}
is an orthogonal system in $A^2(\cW'_\mu )$.
\end{cor}

Thus, we are led to consider the following problem.  
We set $\Delta=\{\z:\, |\z-1|<1\}$ and consider a set of
functions $\{\z^{\lambda_k}\}$, $k=1,2,\dots$. We call the {\em\MS\
  problem for the Bergman space} the question of determining 
necessary and sufficient condition for such a set to
be a {\it complete set} in $A^2(\Delta)$, that is, its linear span to be dense
in $A^2(\Delta)$.  
 The following is Theorem 3.1 in \cite{KPS2}, that gives a sufficient
 condition for the solution of the \MS\ problem for the Bergman space.
\begin{thm}\label{M-S}	
For $k\in\bbN$, $0<a<1$, $c_0>-1$ and $b\in\bbR$, let $\lambda_k=ak+c_0+ib$.
Then $\{ \zeta^{\lambda_k}\}$ is a  complete set in  $A^2(\Delta)$. 
\end{thm}

As a consequence we obtain the following density
result in $A^2(\cW'_\mu)$, which is Theorem 3.1
in \cite{KPS2}.
\begin{thm}\label{complete-set}	 
 Let $\mu>\pi/2$.  
Let $ H _{\ell,j}(z_1,z_2)$ be as in \eqref{H-ell-j-def}. 
Then
$\{ H _{\ell,j}\}_{j\in\bbZ,\,\ell\in\bbN}$, 
is a complete set in $A^2(\cW_\mu')$.
\end{thm}

Notice that the set $\{ H _{\ell,j}:\, j\in\bbZ,\, \ell\in\bbN\}$ is the
union of the two sets in \eqref{two-sets}.  However, such that set is
not an orthogonal set, and we cannot compute the Bergman kernel from
such complete set.

We now divert a bit from our main course to discuss the question of
solving the \MS\ problem. This was done in \cite{PS1,PS2}, however
without finding a complete solution.    In \cite{PS2}
it is proved that the \MS\
problem for the Bergman space is equivalent to characterizing the sets
of
uniqueness of the Hilbert space of holomorphic functions
$\cM^2_\omega(\cR)$ which is the space of holomorphic functions on the
right half-plane $\cR$ such that:
\begin{itemize}
\item[(H)]
$f\in H^2(S_b)$ for every $0<b<\infty$;\smallskip
\item[(B)] $f\in L^2(\ov\cR,d\omega)$;
\end{itemize}
where $H^2(S_b)$ denotes the standard Hardy space on the vertical
strip $\{z\in\bbC:\,  0<\Im z<b\}$, and $\omega $ is the measure on
$\ov\cR$
$$
\omega =\sum_{n=0}^{+\infty} \frac{2^n}{n!}\delta_{\frac{n}{2}}(x)\otimes
dy \,. 
$$ 
Observe thag $\omega$ is a translation invariant measure in
$\ov\cR$. A quite interesting fact is that
such space $\cM^2_\omega(\cR)$ is closely related to a space of
holomorphic functions descovered by T.
Kriete  and D. Trutt, \cite{KT2,KT1}.  The following are the main
results in \cite{PS2} on this problem.

\begin{thm}\label{zero-set-thm}
Let $\{z_j\}\subseteq\cR$, $1\le|z_j|\to+\infty$.  The following properties hold.
\begin{itemize}
\item[(i)] 
If $\{z_j\}$ has exponent of convergence $1$
 and  
upper density  $d^+<\frac12$, then
$\{z_j\}$ is a zero-set for $\cM^2_\omega(\cR)\cap\Hol(\ov\cR)$.\smallskip
\item[(ii)]
If $\{z_j\}$ is a zero-set for $\cM^2_\omega(\cR)\cap\Hol(\ov\cR)$, then
\begin{equation}\label{our-Carleman-cond}
\limsup_{R\to+\infty} \frac{1}{\log R} \sum_{|z_j|\le R}
\Re\big(1/z_j\big) \le 
   \frac2\pi \,.
\end{equation}
\end{itemize}
\end{thm}

We observe that part (ii) in the above theorem follows 
from a generalization of the classical Carleman's formula in the right
half-plane.

\begin{thm}\label{MS-thm}
A sequence $\{z_j\}$ of points in $\cR$ such that $\Re z_j\ge \eps_0$,
for some $\eps_0>0$ and 
that violates condition \eqref{our-Carleman-cond}, is a set of
uniqueness for $\cM^2_\omega(\cR)$.

As a consequence, if $\{z_j\}$ is a sequence as above, the set of
powers $\{ \z^{z_j-1}\}$ is a complete set in $A^2(\Delta)$.  
\end{thm}

\ms

\section*{Final Remarks}

It is worth mentioning that in 
\cite{KPS1} the authors considered the unbounded worm domain
$\cW_\infty$.  In \cite{KPS1} it is proved that the Bergman space
$A^2(\cW_\infty)$ is non-trivial and the Bergman projection is
unbounded on $W^{s,p}(\cW_\infty)$ for all $p\neq2$ and $s>0$.

Moreover, in \cite{BDP} yet another interesting point of view of the
pathological behavior of the worm domain
is considered, in connection with the theory of space–time
singularities associated to the Fefferman metric.\ms

Many questions remain unswered and thus analysis on worm domains is,
and we believe it will remain, a very active area of research.  It
touches function theory and geometry of domains in several complex
variables, holomorphic function spaces, distribution of zeros of
entire functions, regolarity of integral operators, hypoellipticity of
partial differential operators, to name the most significant.

\bibliography{bibWormPlane-1}
\bibliographystyle{amsalpha}

\end{document}